\documentclass[peerreview, 11pt]{IEEEtran}

\usepackage[margin=1in]{geometry}
\usepackage[utf8]{inputenc}
\usepackage{amsmath,bm,amsthm,amssymb}
\usepackage{mathtools}
\usepackage{bbm}
\allowdisplaybreaks
\usepackage{graphicx}
\usepackage[dvipsnames]{xcolor}
\usepackage[english]{babel}
\usepackage[font=small]{caption}
\usepackage{subcaption}
\usepackage{wrapfig}
\usepackage{tabularx}
\usepackage{multirow}
\usepackage{multicol}
\usepackage{thmtools, thm-restate}
\usepackage[normalem]{ulem}
\usepackage[noadjust]{cite}
\usepackage{makecell}
\usepackage{cancel}
\usepackage[colorlinks=true, pdfstartview=FitH, linkcolor=blue, citecolor=blue, urlcolor=blue]{hyperref} 
\usepackage{hyphenat}
\usepackage[linesnumbered,ruled,vlined]{algorithm2e}
\SetKwInput{KwInput}{Input}                
\SetKwInput{KwOutput}{Output}              
\usepackage{algorithmic}

\usepackage[final]{changes}

\newcommand{\floor}[1]{\left\lfloor#1\right\rfloor}

\newtheorem{theorem}{Theorem$\!$}

\newtheorem{lemma}{Lemma$\!$}

\newtheorem{proposition}{Proposition$\!$}

\theoremstyle{definition}
\newtheorem{definition}{Definition$\!$}
\newtheorem{remark}{Remark$\!$}
\usepackage{graphicx}

\usepackage{xcolor}

\newcommand{\id}{\mathrm{id}}

\newcommand{\cB}{\mathcal{B}}
\newcommand{\cC}{\mathcal{C}}

\newcommand{\cG}{\mathcal{G}}
\newcommand{\cH}{\mathcal{H}}

\newcommand{\cK}{\mathcal{K}}

\newcommand{\cS}{\mathbb{S}}

\newcommand{\bbN}{\mathbb{N}}

\newcommand{\mybold}[1]{\bm{#1}}

\newcommand{\bpi}		{\mybold{\pi}}

\newcommand{\bsigma}	{\mybold{\sigma}}
\newcommand{\btau}		{\mybold{\tau}}

\IEEEoverridecommandlockouts

\begin{document}


\title{Improving the Gilbert-Varshamov bound for permutation Codes in the Cayley metric and Kendall $\tau$-Metric}

\author{  
\textbf{The Nguyen}

\IEEEauthorblockA{Department of Mathematics, University of Illinois Urbana-Champaign, Urbana, IL 61801, United States}

{Email: \, thevn2@illinois.edu} 
}

\maketitle
\setcounter{page}{1}
\thispagestyle{plain}
\pagestyle{plain}

\begin{abstract}
The Cayley distance between two permutations $\bpi, \bsigma \in \cS_n$ is the minimum number of \textit{transpositions} required to obtain the permutation $\bsigma$ from $\bpi$. When we only allow adjacent transpositions, the minimum number of such transpositions to obtain $\bsigma$ from $\bpi$ is referred to the Kendall $\tau$-distance. A set $C$ of permutation words of length $n$ is called a $d$-Cayley permutation code if every pair of distinct permutations in $C$ has Cayley distance at least $d$. A $d$-Kendall permutation code is defined similarly. Let $C(n,d)$ and $K(n,d)$ be the maximum size of a $d$-Cayley and a $d$-Kendall permutation code of length $n$, respectively.  In this paper,  we improve the Gilbert–Varshamov bound asymptotically by a factor $\log(n)$, namely
\[ C(n,d+1) \geq \Omega_d\left(\frac{n!\log n}{n^{2d}}\right) \text{ and } K(n,d+1) \geq \Omega_d\left(\frac{n! \log n}{n^{d}}\right).\]
Our proof is based on graph theory techniques. 


\end{abstract}

\section{Introduction}

A permutation code is a non-empty subset of the symmetry group $\cS_n$ equipped with a distance metric. Permutation codes, originally introduced by Slepian~\cite{slepian1965permutation} in 1965 for transmitting data in the presence of additive Gaussian noises, have been extensively studied for their effectiveness in powerline transmission systems against impulsive noise~\cite{chu2004constructions}, as well as in the development of block ciphers~\cite{de2000application}. Especially in recent years, permutation codes under Kendall’s $\tau$-metric, Ulam metric, and Cayley metric have been extensively studied in cloud storage systems, genome resequencing, and the rank modulation scheme of flash memories~\cite{abdollahi2023new, jiang2009rank,jiang2010correcting,farnoud2013error,barg2010codes,buzaglo2015kendall,wang2017snake,zhang2015snake}. Under these metrics, codes are designed to correct transposition errors or translocation errors.

In this paper, we investigate the problem of determining the maximum size of permutation codes in the Cayley metric and Kendall $\tau$-metric. The Cayley distance between two permutations $\bpi, \bsigma \in \cS_n$ is the minimum number of \textit{transpositions} required to obtain the permutation $\bsigma$ from $\bpi$, where transposition is an exchange of two distinct elements. When we only allow adjacent transpositions, i.e. exchange of two distinct adjacent elements, the minimum number of such transpositions to obtain $\bsigma$ from $\bpi$ is referred to the Kendall $\tau$-distance. Let $C(n,d)$ and $K(n,d)$ be the maximum sizes of a permutation code of length $n$ with minimum Cayley distance and Kendall $\tau$-distance \textit{at least} $d$, respectively. Generally, the problem of estimating the maximum size of a code in any metric for given parameters is very difficult. It is known that $C(n,1) = K(n,1) = n!$ and $K(n,2) = n!/2$. Additionally, it is also known that $K(n, d) = 2$ if $\frac{2}{3}\binom{n}{2} \leq  d \leq \binom{n}{2}$ 
(see~\cite[Theorem 10]{buzaglo2015kendall}). Several other bounds are presented in~\cite{barg2010codes,buzaglo2015kendall,wang2017new, wang2021nonexistence,vijayakumaran2016largest}. We refer interested readers to see~\cite{abdollahi2023new, wang2021nonexistence} and references therein for recent results on the size of permutation codes in the Cayley metric and Kendall $\tau$-metric for small distances.    

In general, the Gilbert–Varshamov (GV) and sphere-packing bounds for permutation codes are well-known and generally outperform other bounds for small values of $d$. The main goal of this paper is to provide an asymptotic improvement of the Gilbert-Varshamov bound on the size of permutation codes in the Cayley distance and Kendall $\tau$-metric. Note that there are several works that improved the GV bound for permutation codes of length $n$ in the Hamming distance $d$~\cite{gao2013improvement,tait2013asymptotic,jin2015construction,wang2017new} and block permutation distance~\cite{xu2019new}. To the best of our knowledge, this paper is the first work to improve the GV bound of permutation code in the Cayley metric and Kendall $\tau$-metric. In a recent paper with Wang, Chee, and Vu \cite{wang2024permutation}, an improvement of GV bound is also achieved for permutation codes in the Ulam distance (also known as \textit{deletion distance}), which is recently considered for $q$-ary codes by Alon et al.~\cite{alon2023logarithmically}. 

The rest of this paper is organized as follows. In Section 
 \ref{sec:notation-results}, we review some basic facts about the Cayley metric and Kendall $\tau$-metric together with the statement of our results. In Section \ref{sec:graph-theory}, we introduce some relevant terminologies and results from graph theory and then establish some properties that allow us to compute the distance between two permutations in the Cayley metric. The proofs of main results are presented in Section \ref{sec:proof}. 

\section{Basic Notation and Statement of Results}\label{sec:notation-results}
Denote the set $\{1,2,\ldots,m\}$ as $[m]$. Let $n$ be a positive integer and $\cS_n$ be the set of all permutations on the set $[n]$. Denote $\bpi=(\pi_1,\pi_2,\dotsc,\pi_{n})\in \cS_n$ as a permutation with length $n$. The symbol $\circ$ denotes the composition of permutations. Specifically, for two permutations $\bsigma$ and $\btau$ , their composition, denoted by $\bsigma \circ \btau$, is the permutation with $\bsigma \circ \btau (i) = \bsigma(\btau(i))$ for all $i \in [n]$. All the permutations under this operation form the noncommutative group $\cS_n$ known as the symmetric group on $[n]$ of size $|\cS_n| = n!$.

Given a permutation \( \bpi = (\pi_1, \pi_2, \ldots, \pi_n) \in \cS_n \), the inverse permutation is denoted as \( \bpi^{-1} = (\pi^{-1}_1, \pi^{-1}_2, \ldots, \pi^{-1}_n) \). Here, \( \pi^{-1}_i \) indicates the position of the element \( i \) in the permutation \( \bpi \). For an integer $x \in [n]$, $\pi^{-1}(x)$ indicates the position of $x$ in permutation $\bpi$. We also use cycle decomposition to represent a permutation in this paper; however, whenever we use this, there is no comma in the representation to distinguish it from the above notation. For example, $\bpi = (1, 3, 2, 5, 6, 4)$ has the cycle representation as $(1)(2\,3)(4\,5\,6)$, or $(2\, 3)(4\, 5 \,6)$ for short by deleting all the cycle of length $1$ in the representation. 

\begin{definition}
    For distinct $i,j\in[n]$, a transposition $\varphi(i,j)$ leads to a new permutation obtained by swapping $\pi_i$ and $\pi_j$ in $\bpi$, i.e,
    \begin{equation*}
        \bpi\varphi(i,j)=(\pi_1,\dotsc,\pi_{i-1},\pi_{j},\pi_{i+1},\dotsc,\pi_{j-1},\pi_{i},\pi_{j+1},\dotsc,\pi_n).
    \end{equation*}
    If $|i-j|=1$, $\varphi(i,j)$ is called the \emph{adjacent transposition}.
\end{definition}
Note that each transposition can be viewed as a $2$-cycle permutation. In this paper, we sometime identify $\varphi(i,j)$ with $\bsigma_{i, j} = (i \, j)$, the $2$-cycle permutation on $[n]$ defined as $\bsigma_{i,j}(i) = j, \bsigma_{i,j}(j) = i$, and $\bsigma_{i,j}(k) = k$ for all $k \neq i, j$.  In this way, we have $\bpi \varphi(i, j) =  \bpi \circ \bsigma_{i,j}$ as a composition of two permutations. 

\begin{definition}\label{def:cayley-distance}
    The Cayley distance between two permutations \( \bpi, \bsigma \in \cS_n \), denoted by \( d_{C}(\bpi,\bsigma) \), is defined as the minimum number of transpositions which are needed to change $\bpi$ to $\bsigma$, i.e,
    \begin{equation*}
        d_{C}(\bpi,\bsigma)=\min\{m:\bpi\varphi_{1}\varphi_{2}\dots \varphi_{m}=\bsigma\},
    \end{equation*}
    where $\varphi_{1},\dotsc,\varphi_{m}$ are transpositions and $\bpi\varphi_{1}\varphi_{2}\dots \varphi_{m} = (\dots((\bpi\varphi_{1})\varphi_{2})\dots \varphi_{m})$. 
\end{definition}

\begin{definition}
    The Kendall $\tau$-distance between two permutations \( \bpi, \bsigma \in \cS_n \), denoted by \( d_K(\bpi,\bsigma) \), is defined as the minimum number of adjacent transpositions which are needed to change $\bpi$ to $\bsigma$, i.e,
    \begin{equation*}
        d_{K}(\bpi,\bsigma)=\min\{m:\bpi\varphi_1\varphi_2\dotsm\varphi_m=\bsigma\},
    \end{equation*}
    where $\varphi_1,\dotsc,\varphi_m$ are adjacent transpositions.
\end{definition}

A \textit{$d$-Cayley permutation code} of length $n$ is a subset $\cC$ of $\cS_n$ such that the Cayley distance of any pair of distinct element in $\cC$ is at least $d$. Let $C(n, d)$ be the maximum size of a $d$-Cayley permutation code of length $n$. Define 
\[ \cB_d(\bpi) = \{ \bsigma : d_C(\bpi, \bsigma) \leq d \},\]
and let $V(n,d)$ be size of $|\cB_d(\bpi)|$ for $\bpi \in \cS_n$. Note that $d_C(\cdot, \cdot)$ is invariante, $V(n,d)$ is well-defined, namely it does not depend on $\bpi$. The following proposition provides us the classical bounds for $C(n,d)$.
\begin{proposition}\label{prop:Cayley-GV-bound}
    Let $n>d \geq 0$ be positive integers, $C(n, d+1)$ is bounded as follows
\begin{equation}
    \frac{n!}{V(n,d)} \leq C(n,d+1) \leq \frac{n!}{V(n,\floor{d/2})}.
\end{equation}
\end{proposition}
For $d$ fixed, it is easy to check that $V(n,d) = \Theta_d\left(n^{2d}\right)$, so 
\begin{equation}
    \Omega_d\left(\frac{n!}{n^{2d}}\right) {} \leq C(n,d+1) \leq O_d\left(\frac{n!}{n^{d}}\right).
\end{equation}

The lower bound is known as the Gilbert-Varshamov (GV) bound and the upper bound follows from the fact that all $\cB_{\floor{(d-1)/2}}(\bpi)$ are disjoint for all $\bpi$ in a $d$-Cayley permutation code. Our first main result is to improve the GV bound by a log factor when the distance $d$ is constant. In particular, we have the following theorem.
\begin{theorem}\label{thm-main}
 For $n > d\geq 0$, then 
 $$C(n,d+1) \geq \Omega_d\left(\frac{n! \log n}{n^{2d}}\right).$$
\end{theorem}

Similarly, let $K(n,d)$ be the maximum size of a $d$-Kendall permutation code that is a subset $\cK$ of $\cS_n$ such that the Kendall $\tau$-distance of any pair of distinct elements in $\cK$ is at least $d$. An analog of Proposition \ref{prop:Cayley-GV-bound} gives 
\begin{equation}
        \Omega_d\left(\frac{n!}{n^{d}}\right) {} \leq K(n,d+1) \leq O_d\left(\frac{n!}{n^{\floor{d/2}}}\right).
\end{equation}
We have the following theorem, which gives us an improvement of a log factor, compared with GV bound, for permutation codes in the Kendall $\tau$-metric. 
\begin{theorem}\label{thm:kendall-tau}
For $n > d \geq 0$,  then 
\[ K(n,d+1) \geq \Omega_d\left(\frac{n! \log n}{n^{d}}\right). \]
\end{theorem}

\section{Graph Theory Perspective}\label{sec:graph-theory}

\subsection{Independence number of a graph with few triangles}
An independent set of a graph $G$ is a subset of vertices such that no two vertices in the subset represent an edge of $G$. For a graph $G$, let $\alpha(G)$ be the independence number of $G$, i.e. the size of a largest independent set in $G$. We recall the following folklore result in graph theory. 
\begin{proposition}\label{prop:indep-number}
Let $G$ be a graph on $n$ vertices with maximum degree $\Delta$. Then $G$ contains an independent set of size $n/(\Delta  + 1)$, in other words $\alpha(G) \geq n/(\Delta +1)$.  
\end{proposition}

When the graph is \textit{locally sparse}. i.e only few triangles, Ajtai, Koml\'os and Szemer\'edi \cite{ajtai1980note} showed that Proposition \ref{prop:indep-number} can be improved by a log factor. More precisely, they showed that for any $\varepsilon > 0$ and any graph $G$ on $N \geq 1$ vertices with average degree $d$ and $T < N d^{2 - \varepsilon}$ triangles, then $\alpha(G) > \Omega_{\varepsilon}\left((N/d)\log d\right)$. This line of work has led to several important developments in extremal graph theory and Ramsey theory. We mention that we can easily replace average degree $d$ by maximum degree $\Delta$ in the Ajtai, Koml\'os, and Szemer\'edi result. In this paper, we will use the following form, which can be found in~\cite{bollobas1998random}. 
\begin{lemma}[\cite{bollobas1998random}, p.296]\label{lemma-graph}
Let $G$ be a graph with maximum degree $\Delta$ ($\Delta \geq 1$) and suppose that $G$ has $T$ triangles. Then 
\begin{equation}
    \alpha(G) \geq \frac{|V(G)|}{10\Delta}\left(\log \Delta - \frac{1}{2}\log\left(\frac{T}{|V(G)|}\right)\right).
\end{equation}
\end{lemma}

\subsection{Cayley graph on $\cS_n$}

In this section, we discuss about the structure of Cayley graph on $\cS_n$ generated by sets of $2$-cycles or transpostions. Given a subset $D$ of a group $(G, \cdot)$, we write $\Gamma_{G}(D)$ to denote the Cayley graph of $G$ with respect to the generator set $D$; that is, the graph on the vertex set $G$ that contains an directed edge $(x, y)$ if and only if $x \cdot y^{-1} \in D$.  
\begin{definition}\label{def:graph-dis}
For $n \in \bbN$ consider the Cayley graph $\Gamma_{\cS_n}(D)$ where $D$ is the set of all transposition permutations, i.e $2$-cycle permutations 
\[ D  = \{ \bsigma_{i,j} \,|\, 1 \leq i, j \leq n\}.\]
Then the \textit{Cayley distance} is the distance on $\cS_n$ which is the \textit{graph distance} of the Cayley graph $\Gamma_{\cS_n}(D)$. 
\end{definition}
Note that $\Gamma_{\cS_n}(D)$ is undirected graph since if $\bpi \circ \sigma^{-1} = \bsigma_{i,j}$ then $\bsigma \circ  \bpi^{-1} = \bsigma_{i,j}$. Clearly, the Definition \ref{def:graph-dis} is equivalent to the Definition \ref{def:cayley-distance}. This perspective allows us to apply the extensive research on Cayley graph to analyze Cayley permutation codes. The following proposition provides a formula to compute the Cayley distance between two permutations, which can be found in \cite[p.118]{diaconis1988group}. We include  the proof here for sake of convenience and completeness.

\begin{proposition}[Cayley's formula]\label{prop:cayley}
    The Cayley distance distance \( d_C(\bpi,\bsigma)\) between \( \bpi \) and \( \bsigma \) equals  $n$ minus the number of permutation cycles in $\bpi \circ \bsigma^{-1}$, i.e 
    \[ d_C(\bpi, \bsigma) = n - \left|\mathrm{Cycles}(\bpi \circ \bsigma^{-1})\right|. \]
     where \( \mathrm{Cycles}(\bpi) \) is the set of permutation cycles in $\bpi$, including $1$-cycles. 
\end{proposition}
\begin{IEEEproof}
Firstly, it follows from the definition that the Cayley metric is a left/right invariant metric, i.e.
\[ d_C(\bpi \circ \bsigma, \btau \circ \bsigma) = d_C(\bpi, \btau) = d_C(\bsigma \circ\bpi , \bsigma\circ\btau)\]
 for all $\bpi, \bsigma, \btau \in \cS_n$. Therefore, it suffices to prove the Proposition \ref{prop:cayley} for the case that $\bsigma = \id$, the identity permutation of $\cS_n$. We now show that for any permutation $\bpi$, the Cayley distance $d_C(\bpi, \id)$ equals $n$ minus the number of cycles of $\bpi$. 
 
 Indeed, we observe that a cycle permutation of $k$ elements (a $k$-cycle) is a composite of $k-1$ transpositions. For example, the $k$-cycle permutation $(1 \, 2\, 3 \, \dots \, k)$ can be written as $\bsigma_{1, 2} \circ \bsigma_{2, 3} \circ \dots \circ \bsigma_{k-1, k}$. We also have  $k-1$ is the minimum number of transpositions such that their composition is a $k$-cycle. This can be easily seen by induction on $k$. Additionally, any permutation $\bpi$ can be expressed as the product of disjoint cycles. All together implies that $$d_C(\bpi, \id) =  n - |\mathrm{Cycles}(\bpi)|.$$
\end{IEEEproof}

\section{Proof of Theorem \ref{thm-main} and Theorem \ref{thm:kendall-tau}} \label{sec:proof}

The proof idea of Theorem \ref{thm-main} and Theorem \ref{thm:kendall-tau} relies on techniques from graph theory by constructing auxiliary graphs on the set of permutations such that each permutation code is an independent set of these graphs. Therefore, studying the size of codes is equivalent to analyzing the independence number of the auxiliary graph. Hence, the Gilbert-Varshamov bound is simply the application of Proposition \ref{prop:indep-number} on this graph. By this terminology, if we can show that the auxiliary graph has few triangles, we have an improvement by a log factor compared with the classical GV bound. This idea was initially used by Jiang and Vardy~\cite{jiang2004asymptotic} to improve GV bound for binary codes with Hamming distance, followed by subsequent papers in other settings, see~\cite{alon2023logarithmically, gao2013improvement,tait2013asymptotic,wang2024permutation} for example. 

\subsection{Proof of Theorem \ref{thm-main}}

We define the graph $\cG_{n, d}$ with vertex set $\cS_n$, and two permutations are connected by an edge if their Cayley distance is most $d$. Therefore, a $(d+1)$-Cayley permutation code is an independent set of $\cG_{n, d}$. Note that $|V(\cG_{n,d})| = n!$ and $\Delta = O_{d}(n^{2d})$. Thus, in order to prove Theorem \ref{thm-main}, by using Lemma \ref{lemma-graph}, it suffices to show that the number of triangles in $\cG_{n, d}$ is $O_d(n! n^{4d - \varepsilon})$ for some $\varepsilon > 0$. We will prove this holds for $\varepsilon = 1$. 

\begin{lemma}\label{lem:num-triangle}
 The number of triples $(\bsigma, \bpi, \btau) \in (\cS_n)^3$ such that $d_C(\bsigma,  \bpi) \leq d, d_C(\bpi, \btau) \leq d$, and $d_C(\bsigma, \btau) \leq d$ is at most $O_d(n! n^{4d-1})$. 
\end{lemma}
\begin{IEEEproof}
For $\bsigma \in S_n$, we will count the number of permutations $\bpi$ and $\btau$ such that $(\bsigma, \bpi, \btau)$ forms a triangle in $\cG_{n, d}$. Note that $(\bsigma, \bpi, \btau)$ is uniquely determined by $\bsigma$ and sequences $S_1 = (\varphi_1,\varphi_2,\dotsc, \varphi_h), S_2 = (\varphi'_1,\varphi'_2, \dotsc \varphi'_\ell)$ of transpositions for which $\bpi = \bsigma \varphi_1\varphi_2\dotsm\varphi_h$ and $\btau = \bpi\varphi'_1\varphi'_2\dotsm\varphi'_\ell$. For the sake of brevity, we denote $\bpi\circ S=\bpi\varphi_1\varphi_2\dotsm\varphi_m$, where the transposition sequence is $S=(\varphi_1,\varphi_2,\dotsc,\varphi_m)$.

Since $d_C(\bsigma, \bpi) \leq d$ and $d_C(\bpi, \btau) \leq d$, we can choose $S_1, S_2$ such that $|S_1|, |S_2| \leq d$. Hence, we can combine $S_1$ and $S_2$ to obtain a sequence $S$ of transpositions with size $|S_1| + |S_2| \leq 2d$ such that $\btau = \bsigma \circ S$. Now, for a such sequence $S$, there are at most $O_d(1)$ possibilities for $S_1$ and $S_2$ that produce $S$. Therefore, for a given $\bsigma$, the number of triangles $(\bsigma, \bpi, \btau)$ is at most $O_d(1)$ times the number of ways to choose a transposition sequence $S$ with $|S| \leq 2d$ such that $\btau = \bsigma \circ S$ and $d_C(\bsigma, \btau) \leq d$. Hence, it suffices to show that for a given permutation $\bsigma$, the number of transposition sequences $S$ with $|S| \leq 2d$ such that $\btau = \bsigma \circ S$ and $d_C(\bsigma, \btau) \leq d$ is at most $O_d(n^{4d-1})$. 

By symmetry, we can assume that $\bsigma = \mathrm{id}$ and $|S| = 2d$, here $\mathrm{id}$ stands for identity permutation, i.e $\mathrm{id}= (1, 2, \dots, n)$. Let $S =  \varphi(a_1, b_1) \circ \varphi(a_2, b_2)  \circ \dots \circ \varphi(a_{2d}, b_{2d})$ be a transposition sequence such that $d_C(\mathrm{id}, \mathrm{id} \circ  S)  \leq d$. If $a_i,b_j$ are distinct for all $1 \leq i, j \leq d$, then $\btau = \mathrm{id} \circ S$ is determined by $\btau(a_i) = b_i$ and $\btau(b_i) = a_i$ for all $1 \leq i \leq 2d$ and $\btau(x) = x$ for all other positions, i.e $\btau = (a_1 \, b_1) (a_2 \, b_2) \dots (a_{2d} \, b_{2d})$, here we write a permutation as a product of cycles. This is because each transposition does not affect the positions of other elements and $a_i, b_j$ are distinct for all $1\leq i, j \leq d$. It follows that $\btau$ has $n-2d$ disjoint cycles, so $d_C(\mathrm{id}, \btau) = 2d$ by Proposition \ref{prop:cayley}, a contradiction. Therefore, the multiset $\{a_i, b_j : 1 \leq i, j \leq 2d\}$ has at most $4d-1$ distinct values. Hence, the number of choices for $S$ is at most 
\[ \binom{n}{4d-1} (4d-1)^{4d} = O_d(n^{4d-1}).\]
This completes the proof of Lemma \ref{lem:num-triangle}.
\end{IEEEproof}

\subsection{Proof of Theorem \ref{thm:kendall-tau}}
Let $\cH_{n,d}$ be the graph with $V(\cH_{n,d}) = \cS_n$ and there is an edge between two permutation $\bsigma$ and $\btau$ if $d_K(\bsigma,  \btau) \leq d$.  Hence, a $(d+1)$-Kendall permutation code is an independent set of $\cH_{n,d}$. It is easy to see that $|V(\cH_{n,d})| = n!$ and $\Delta = O_d(n^d)$. Thus, Theorem \ref{thm:kendall-tau} is followed by Lemma \ref{lemma-graph} and the following lemma, which is an analog of Lemma \ref{lem:num-triangle} for Kendall $\tau$-distance.  

\begin{lemma}\label{lem:num-triangle-kendall}
 The number of triples $(\bsigma, \bpi, \btau) \in (\cS_n)^3$ such that $d_K(\bsigma,  \bpi) \leq d, d_K(\bpi, \btau) \leq d$, and $d_K(\bsigma, \btau) \leq d$ is at most $O_d(n! n^{2d-1})$. 
\end{lemma}

\begin{IEEEproof}
The proof is almost identical to Lemma \ref{lem:num-triangle}, in which we consider adjacent transpositions instead of general transpositions. It suffices to prove for a given permutation $\bsigma$, the number of sequences $S$ of adjacent transpositions with $|S| \leq 2d$ such that $\btau = \bsigma \circ S$ and $d_K(\bsigma, \btau) \leq d$ is at most $O_d(n^{2d-1})$. 

Without loss of generality, we can assume $\bsigma = \mathrm{id}$ and $|S| = 2d$. Let $S = \varphi(i_1, i_1+1)  \circ \varphi(i_2, i_2+1) \circ \dots \circ \varphi(i_{2d}, i_{2d} +1)$ be a sequence of adjacent transpositions such that $d_K(\mathrm{id}, \id \circ S) \leq d$. If $|i_a -  i_{b}| \geq 2$ for all $1 \leq a < b \leq 2d$ then we have 
\[ \btau = \mathrm{id} \circ S = (i_1 \, i_1+1) (i_2\, i_2+1) \dots (i_{2d}\, i_{2d}+1).\]
However, in this case, we have $d_K(\mathrm{id}, \btau) = 2d$, a contradiction. Therefore, there exists $1 \leq a < b \leq 2d$ such that $|i_a - i_b| \leq 1$, thus in $S$ we should have either $\{(i_a , i_a +1), (i_a+1 , i_a +2)\}$, $\{(i_a , i_a +1), (i_a-1 , i_a)\}$ or $(i_a \,i_a +1)$ twice. Therefore, we have 
\[ |S| \leq n \binom{n}{2d-2} O_d(1) = O_d(n^{2d-1}),\]
which completes the proof of Lemma \ref{lem:num-triangle-kendall}.
\end{IEEEproof}

\begin{remark}
The exponent $4d-1$ and $2d-1$ in Lemma \ref{lem:num-triangle} and Lemma \ref{lem:num-triangle-kendall} can be improved by a more precise argument. However, we present the above proofs for simplicity since this is enough to obtain our main reults. 
\end{remark}

\section*{Acknowledgement}
The author would like to thank Jozsef Balogh,  Bob Krueger, Quy Dang Ngo, Van Khu Vu, Micheal Wigal, and Ethan White for useful discussions. The author is grateful to Ethan White for a careful reading of the manuscript. The author was partially supported by a David G. Bourgin Mathematics Fellowship. 

\bibliographystyle{IEEEtran}
\bibliography{references.bib}
\end{document}